\documentclass[10pt,twoside,a4paper]{amsart}
\usepackage{amsmath}
\usepackage{a4}

\usepackage{amssymb}
\usepackage{setspace}
\newcommand{\WFN}{\operatorname{WFN}}
\newcommand{\frow}{{}^\frown }

\newcommand{\supt}{\operatorname{supt}}

\newcommand{\trcl}{\operatorname{trcl}}

\newcommand{\sk}{\operatorname{sk}}

\newcommand{\cf}{\operatorname{cf}}

\newcommand{\Fn}{\operatorname{Fn}}

\newcommand{\dom}{\operatorname{dom}}

\newcommand{\cov}{\operatorname{cov}}

\newcommand{\id}{\operatorname{id}}

\newcommand{\non}{\operatorname{non}}

\newcommand{\Mgood}[1]{\mathcal M_{#1}}

\newcommand{\forces}{\Vdash}
\newcommand{\SEP}{\operatorname{SEP}}
\newcommand{\IDP}{\operatorname{IDP}}

\newcommand{\card}[1]{\mid\!\!{#1}\!\!\mid}
\newcommand{\cnt}[1]{[{#1}]^{\leq{\aleph_0}}}
\newtheorem{thm}{Theorem}[section]
\newtheorem{corollary}[thm]{Corollary}
\newtheorem{lemma}[thm]{Lemma}

\theoremstyle{definition}
\newtheorem{defn}[thm]{Definition}

\newcommand{\setof}[2]{\{#1\,:\,#2\}}
\newcommand{\psof}[1]{{\mathcal P}\/(#1)}
\newcommand{\pairof}[1]{(#1)}
\newcommand{\seqof}[2]{(#1)_{#2}}
\newcommand{\cardof}[1]{\mathopen{|\,}#1\mathclose{\,|}}
\newcommand{\st}{such that}

\newcommand{\calC}{{\mathcal C}}
\newcommand{\calH}{{\mathcal H}}
\newcommand{\calM}{{\mathcal M}}

\begin{document}

\title{Remarks on a paper by Juh\'asz and Kunen}
\author{Saka\'e Fuchino}
\address[Fuchino]{Dept.~of Natural Science and Mathematics\\
  College of Engineering\\
  Chubu University\\
  Kasugai Aichi 487-8501 JAPAN}
\email{fuchino@isc.chubu.ac.jp}
\author{Stefan Geschke}
\address[Geschke]{II.~Mathematisches Institut\\Freie Universit\"at
Berlin\\Arnimallee 3\\14195 Berlin\\Germany}
\email{geschke@math.fu-berlin.de}
\date{September 12, 2003}
\thanks{The first author thanks Hiroshi Sakai for pointing out a gap in the proof of 
an earlier version of Theorem 8.14} 

\begin{abstract}We give an equivalent, but simpler formulation of the
axiom SEP introduced by Juh\'asz and Kunen in \cite{juku}.  Our
formulation shows that many of the consequences of the weak
Freese-Nation Property of $\mathcal P(\omega)$ studied in
\cite{fugeso} already follow from SEP.
We show that it is consistent that
SEP holds while $\mathcal P(\omega)$ fails to have the
$(\aleph_1,\aleph_0)$-ideal property introduced in \cite{dow-hart}.
This answers a question addressed independently by
Fuchino and by Kunen.
We also consider some natural variants of SEP and show that certain changes
in the definition of SEP do not lead to a different principle.  
This answers a question addressed by Blass.  
\end{abstract}

\maketitle

\section{Introduction}
In \cite{juku}, Juh\'asz and Kunen introduced a new principle, $\SEP$,
which is designed to capture some consequences of CH and holds in Cohen
models, i.e., models of set theory obtained by adding Cohen reals to
models of CH.
The principle
says that for sufficiently large $\chi$ there are
many good elementary submodels $M$ of the structure $(\mathcal
H_\chi,\in)$ such that $\mathcal P(\omega)\cap M$ is nicely embedded
in $\mathcal P(\omega)$.

Different notions of {\em good} and {\em many} elementary submodels and
{\em nice} embeddings lead to different principles.  In this paper we
study the relation between $\SEP$ and the
$(\aleph_1,\aleph_0)$-ideal
property of Dow and Hart \cite{dow-hart}, which is defined by the same
pattern.  The corresponding
embeddings are $sep$-embeddings and $\sigma$-embeddings.

\begin{defn} Let $A$ and $B$ be Boolean algebras such that $A\leq B$.
Then $A\leq_{sep}B$ if and only if for every $b\in B$ and every uncountable set
$T\subseteq A\restriction b$ there is $a\in A\restriction b$ such that
$\{c\in T:c\leq a\}$ is uncountable.
Let $A\leq_{\sigma}B$ if and only if for every $b\in B$, $A\restriction b$ has a
countable cofinal subset.
\end{defn}

Observe that $A\leq_\sigma B$ implies $A\leq_{sep} B$.

The good elementary submodels will be the same for both $\SEP$ and the
$(\aleph_1,\aleph_0)$-ideal property:
later we show that $\leq_{sep}$ can
be replaced by $\leq_\sigma$ in the following definition of $\SEP$.
The real
difference between $\SEP$ and the $(\aleph_1,\aleph_0)$-ideal property
lies in the interpretation of {\em many}.

Extending the notions defined by Juh\'asz and Kunen, we regard SEP and
the $(\aleph_1,\aleph_0)$-ideal property as properties of general Boolean
algebras, not only of $\mathcal P(\omega)$.

\begin{defn}
For a cardinal $\chi$ let $\Mgood{\chi}$ be the set of elementary
submodels $M$ of $\mathcal H_\chi$ such that $\card{M}=\aleph_1$ and
$[M]^{\aleph_0}\cap M$ is cofinal in $[M]^{\aleph_0}$.
For a Boolean algebra $A$, $\SEP(A)$ is the statement ``for all
sufficiently large regular $\chi$ there are cofinally many
$M\in\Mgood{\chi}$ such that $A\cap M\leq_{sep}A$''.
The axiom $\SEP$ introduced in \cite{juku}
is $\SEP(\mathcal P(\omega))$.

$A$ has the $(\aleph_1,\aleph_0)$-ideal property (IDP) if and only if for all
sufficiently
large regular $\chi$ and all $M\in\Mgood{\chi}$ with $A\in M$,
$A\cap M\leq_\sigma A$ holds.  We write $\IDP(A)$ if $A$ has the IDP.
\end{defn}

As usual, we identify an algebraic structure
$(A,f_1,\dots,f_n,R_1,\dots,R_m)$ with its underlying set $A$.
It should be made clear that by ``$A\in M$''
we really mean $(A,f_1,\dots,f_n,R_1,\dots,R_m)\in M$.
The role of the cardinal $\chi$ is analyzed in Section \ref{section-eight}, where 
we calculate precisely when a regular cardinal is ``sufficiently large'' in the
definition of SEP. 

Since $\leq_\sigma$ is stronger than $\leq_{sep}$, for a Boolean algebra
$A$, $\IDP(A)$ implies $\SEP(A)$.  
In Section 2 we show that the relation $\leq_{sep}$ in the definition of $\SEP$ 
can be replaced by $\leq_\sigma$, i.e., SEP and IDP are very similar.  
In Section \ref{section-eight} we observe that it does not make a difference in the 
definition of SEP if we replace ``there are cofinally many $M\in\Mgood{\chi}$''
by ``there is $M\in\Mgood{\chi}$'' or by ``there are stationarily many (in $[\mathcal H_\chi]^{\aleph_1}$) $M\in\Mgood{\chi}$''.  
This shows the importance of the results of Section 6 and Section 7,
namely that SEP is really weaker than IDP.

If the universe is not 
very complex, that is, if some very weak version of the
$\square$-principle together with $\cf([\mu]^{\aleph_0})=\mu^+$ holds for
all singular cardinals $\mu$ of countable cofinality, or if the Boolean
algebras under consideration are small,
then the $(\aleph_1,\aleph_0)$-ideal property
is equivalent to the weak Freese-Nation property studied in
\cite{fukoshe}: 

\begin{defn} A Boolean algebra
$A$ has the weak Freese-Nation property
($\WFN$) if and only if there is a function $f:A\to[A]^{\leq\aleph_0}$ such that
for all $a,b\in A$ with $a\leq b$ there is $c\in f(a)\cap f(b)$ with
$a\leq c\leq b$.  $f$ is called a WFN-function for $A$.
We write $\WFN(A)$ for ``$A$ has the WFN''.
\end{defn}

It is easy to check that
if $f$ is a WFN-function for $A$ and $B\leq A$ is closed under $f$, then
$B\leq_\sigma A$.
In \cite{fukoshe} Fuchino, Koppelberg, and Shelah characterized the WFN
using elementary submodels and $\sigma$-embeddings.  They showed

\begin{thm}\label{charwfn}
A Boolean algebra $A$ has the WFN if and only if for all sufficiently large $\chi$
and all $M\preccurlyeq\mathcal H_\chi$ with $\aleph_1\subseteq M$ and
$A\in M$, $A\cap M\leq_\sigma A$.
\end{thm}

Since $\aleph_1\subseteq M$ for all $M\in\Mgood{\chi}$, it is clear that
for a Boolean algebra $A$, $\IDP(A)$ follows
from $\WFN(A)$.
As mentioned above, IDP and WFN
are equivalent in many cases.
Some of these cases are captured by the following lemma, which
follows from the results in \cite{fuso}.

\begin{lemma}\label{wfn} If $A$ is a Boolean algebra of size
$<\aleph_\omega$ or if $0^\sharp$ does not exist, then
$\IDP(A)$ holds if and only if $\WFN(A)$ does.
\end{lemma}

The formulation using $0^\sharp$ is chosen here
just for simplicity.
As mentioned above, what is really needed is only a certain very weak assumption
at singular cardinals of countable cofinality. On the other hand, it is
known that the lemma does not hold without any such additional
assumption (see \cite{fuso} or \cite{fugesheso}).

In \cite{fugeso} many interesting consequences of $\WFN(\mathcal
P(\omega))$ have been found.  Concerning the combinatorics of the reals,
a universe satisfying $\WFN(\mathcal P(\omega))$ behaves very similar to a
Cohen model.
In particular, the values of the popular cardinal invariants of the
continuum, that is, those studied in \cite{blass}, have the same values
in a model with
$\WFN(\mathcal P(\omega))$ as in a Cohen model with the same size of the
continuum.

Our characterization of $\SEP$ in terms of $\leq_\sigma$ rather than $\leq_{sep}$ 
shows that the axiom $\SEP(\mathcal P(\omega))$ is sufficient to determine at
least some of the smaller cardinal invariants of the continuum.

Juh\'asz and Kunen already proved that another consequence of
$\WFN(\mathcal P(\omega))$, the principle $C_2^s(\omega_2)$ introduced in
\cite{jusosze}, follows from the weaker assumption $\SEP(\mathcal
P(\omega))$.

\section{$\SEP$ is similar to $\IDP$}

As mentioned above, $\leq_\sigma$ implies $\leq_{sep}$.  If the
subalgebras under consideration are of size $\leq\aleph_1$, then the two
relations are in fact the same.

\begin{lemma}  Let $A$ and $B$ be Boolean algebras with $A\leq_{sep} B$
and $\card{A}=\aleph_1$.  Then $A\leq_\sigma B$.
\end{lemma}

\begin{proof}
Let $b\in B$ and assume for a contradiction that $A\restriction b$ is not
countably generated.  Let $(a_\alpha)_{\alpha<\omega_1}$ enumerate
$A\restriction b$.  By recursion on $\alpha<\omega_1$ define a sequence
$(c_\alpha)_{\alpha<\omega_1}$ in $A\restriction b$ such that for all
$\alpha<\omega_1$, $a_\alpha\leq c_\alpha$ and
$c_\alpha$ is not in the ideal of $A$ generated by
$\{c_\beta:\beta<\alpha\}$.

Now let $T=\{c_\alpha:\alpha<\omega_1\}$.  We claim that $T$ is a
counterexample to $A\leq_{sep} B$.  For let $a\in T$.  Then there is
$\beta<\omega_1$ with $a=a_\beta$.  Thus $a_\beta\leq c_\beta$.  By
the construction of the sequence $(c_\alpha)_{\alpha<\omega_1}$, there are
only countably many elements of $T$ below $c_\beta$.  Hence, there are
only countably many elements of $T$ below $a$, contradicting
$A\leq_{sep}B$.
\end{proof}

Thus, we have 

\begin{corollary}\label{sepvsidp} For every Boolean algebra $A$, $\SEP(A)$
holds if and only if
for all sufficiently large regular $\chi$ there are cofinally many $M\in
\Mgood{\chi}$ with $A\cap M\leq_\sigma A$.
\end{corollary}

In this characterization, it is easily seen that many interesting
consequences of $\WFN(\mathcal P(\omega))$ already follow from
$\SEP(\mathcal P(\omega))$.  In the proofs of most of the results in
\cite{fugeso} it is only used that under $\WFN(\mathcal P(\omega))$,
for some sufficiently large $\chi$ there
are cofinally many $M\in\Mgood{\chi}$ with $\mathcal P(\omega)\cap
M\leq_\sigma\mathcal P(\omega)$.  The following theorem collects
some of the consequences of $\SEP(\mathcal P(\omega))$ that follow from
the arguments given in \cite{fugeso}.

A subset of $[\omega]^{\aleph_0}$ is called groupwise dense if
it is closed under taking almost subsets and
non-meager with respect to the topology on
$[\omega]^{\aleph_0}$ inherited from $2^\omega$
when identifying $[\omega]^{\aleph_0}$ with a subset of $2^\omega$.

\begin{thm} Assume $\SEP(\mathcal P(\omega))$.  Then the following
cardinal
invariants of the continuum are $\aleph_1$.
\begin{enumerate}\item[1.] $\non(\mathcal M)$, the smallest size of a
non-meager subset of $\mathbb R$ and
\item[2.] $\mathfrak g$, the smallest size of a family of groupwise dense
subsets of $[\omega]^{\aleph_0}$ with empty intersection.
\end{enumerate}
Moreover, if CH fails, then $\cov(\mathcal M)$, the minimal size of a
family of non-meager subsets of $\mathbb R$ covering $\mathbb R$, is
at least $\aleph_2$.
\end{thm}

It was also proved in \cite{fugeso} that $\WFN(\mathcal P(\omega))$
implies that $\mathfrak a$, the smallest size of a maximal almost disjoint
family in $\mathcal P(\omega)$, is $\aleph_1$.  In the proof it is
sufficient to assume $\IDP(\mathcal P(\omega))$.
The situation with $\mathfrak a$ under $\SEP(\mathcal P(\omega))$ is
more subtle and
will be discussed in section \ref{section-eight}.

For more information about these cardinal invariants see e.g.\
\cite{blass}.

\section{$\SEP(\mathcal P(\omega))$ holds in Cohen models}

In \cite{fugesheso} it was shown that $\WFN(\mathcal P(\omega))$ can fail
in a Cohen model, assuming the consistency of some very large cardinal.
In \cite{fuso} it was shown that large cardinal assumptions are
necessary for this.

$\SEP(\mathcal P(\omega))$ however,
is always true in a Cohen model.  We include a proof of this fact
(Theorem \ref{sepcohen}).
It follows that $\WFN(\mathcal P(\omega))$ does not
follow from $\SEP(\mathcal P(\omega))$ in ZFC (assuming the consistency of certain large cardinals).
We do not know whether $\IDP(\mathcal P(\omega))$ is always true in a
Cohen model.

In section \ref{section-seven}
we shall show, without large cardinal assumptions,
that $\SEP(\mathcal P(\omega))$ does not imply
$\IDP(\mathcal P(\omega))$.

\begin{thm}\label{sepcohen}
Let $V$ be a model of CH and suppose that
$G$ is
$\Fn(\kappa,2)$-generic over $V$.
Then $V[G]\models\SEP(\mathcal P(\omega))$.
\end{thm}

The proof of this theorem relies on the following series of lemmas.  The
first lemma was
proved in \cite{shelah-steprans}.

\begin{lemma}\label{cohenreal} Let $M$ be a transitive model of set theory
and let
$x\subseteq\omega$ be a Cohen real over $M$.  Then
${\mathcal P}(\omega)\cap M\leq_{\sigma}{\mathcal P}(\omega)\cap M[x]$.
\end{lemma}

\begin{lemma}\label{submodel} Let $P$ be a c.c.c.~partial order and let
$M\in\Mgood{\chi}$ be such that $P\in M$.  Let $G$ be $P$-generic
over the ground model $V$.  Then $M[G]\in\Mgood{\chi}^{V[G]}$.
\end{lemma}

\begin{proof} Let $M$ and $G$ be as above.  Then
$M[G]\preccurlyeq\mathcal H_\chi[G]$ by c.c.c.~of $P$.
Let $X\subseteq M[G]$ be countable.
In $V[G]$, there is a countable set $C\subseteq M$ of $P$-names such that
$X\subseteq\{\dot{x}_G:\dot{x}\in C\}$.  Again by the c.c.c.~of $P$,
we may assume that $C\in V$ and $C$ is countable in $V$.  By
$M\in\Mgood{\chi}^V$,
we may assume
$C\in M$.  Since $M$ contains a name
for $G$,  $G\in M[G]$ and thus
$X\subseteq\{\dot{x}_G:\dot{x}\in C\}\in M[G]$.
\end{proof}

\begin{lemma}
\label{CH-M}
Assume {\rm CH}. Then $M\in\calM_\chi$ implies that
$[M]^{\leq\aleph_0}\subseteq M$.

In particular, for all $M\in\calM_\chi$, if $P\subseteq M$ 
is a c.c.c.~partial order, then every nice $P$-name for a subset of $\omega$
is contained in $M$.
\end{lemma}
\begin{proof}
Suppose that $x\in[M]^{\aleph_0}$. By $M\in\calM_\chi$ there is
$y\in[M]^{\aleph_0}\cap M$ such that $x\subseteq y$. By CH there is a
surjection $f:\omega_1\rightarrow\psof{y}$. By elementarity, there is
such $f$ in $M$. Let $\alpha<\omega_1$ be such that $f(\alpha)=x$.
Since $\omega_1\subseteq M$, we have $\alpha\in M$ and thus
$x=f(\alpha)\in M$.

Now let $P\subseteq M$ be a c.c.c.~partial order.
Since $P$ is c.c.c., 
every nice $P$-name for a subset of $\omega$ is 
a countable subset of $\{\check n:n\in\omega\}\times P\subseteq M$.
By the first part of the lemma, every such set is an element of $M$.
\end{proof}

\begin{lemma}\label{submodelPomega}
Let $V$ and $G$ be as in Theorem
\ref{sepcohen}. Let
$\chi\in V$ be large enough and let $M\in(\Mgood{\chi})^V$ be such that
$\kappa\in M$.
Then $\mathcal P(\omega)\cap M[G]=\mathcal P(\omega)\cap
V[G\cap\Fn(\kappa\cap M,2)]$.
\end{lemma}

\begin{proof} First let $x\in\mathcal P(\omega)\cap M[G]$.
In $M$ there is a name $\dot{x}$
for $x$ which is a nice name for a subset of $\omega$.
Since $\Fn(\kappa,2)$ satisfies the c.c.c., in $V$ there is a countable set
$X\subseteq\kappa$ such that $\dot{x}$ is an $\Fn(X,2)$-name.
We can find such an $X$ in $M$.  Since $\aleph_0\subseteq M$,
$X\subseteq M$.  
Therefore $\dot{x}$ is an $\Fn(\kappa\cap M,2)$-name and
thus $x\in V[G\cap\Fn(\kappa\cap M,2)]$.

For the converse let $x\in V[G\cap\Fn(\kappa\cap M,2)]$.
Pick a nice $\Fn(\kappa\cap M,2)$-name $\dot{x}$ for $x$.
Clearly, $\Fn(\kappa\cap M,2)\subseteq M$.  
Therefore Lemma \ref{CH-M} applies, and we get $\dot{x}\in M$.
This shows $x\in M[G]$.
\end{proof}

\begin{proof}[Proof of Theorem \ref{sepcohen}]
We argue in $V[G]$.  Let $\chi$ be sufficiently large and let
$M\in\Mgood{\chi}$.  It is sufficient to show that there is
$M^\prime\in\Mgood{\chi}$ with $M\subseteq M^\prime$ and $\mathcal
P(\omega)\cap
M^\prime\leq_\sigma\mathcal P(\omega)$.

Since $\Fn(\kappa,2)$ has the c.c.c., in $V$ there is a set $X$ of size $\aleph_1$ of
$\Fn(\kappa,2)$-names  such that every element of $M$
has a name in $X$.  Let $N\in\Mgood{\chi}^V$ be such that $X\subseteq N$.
Clearly, $M\subseteq N[G]$.  By Lemma \ref{submodel},
$N[G]\in\Mgood{\chi}$.  By Lemma \ref{submodelPomega}, $\mathcal
P(\omega)\cap N[G]=\mathcal P(\omega)\cap V[G\cap\Fn(\kappa\cap N,2)]$.
Since $V[G]$ is an $\Fn(\kappa\setminus N,2)$-generic extension over
$V[G\restriction\Fn(\kappa\cap N,2)]$, it follows from Lemma
\ref{cohenreal} that $\mathcal P(\omega)\cap
N[G]\leq_\sigma\mathcal P(\omega)$.  This shows that, in $V[G]$, the set of
$M'\in{\mathcal M}_\chi$ with
${\mathcal P}(\omega)\cap M'\leq_\sigma{\mathcal P}(\omega)$ is cofinal
in $[{\mathcal H}_\chi]^{\aleph_1}$.
\end{proof}

\section{$\IDP$ for partial orders}
In \cite{fukoshe} the $\WFN$ has been defined for partial orders, not only
for Boolean algebras.  In this section and the next, we do the same
for $\SEP$ and $\IDP$.
We have to liberalize our definition of $\leq_\sigma$.

\begin{defn} Let $P$ and $Q$ be partial orders with $Q\leq P$, i.e.,
$Q\subseteq P$ and the orders on $Q$ and $P$ agree on $Q$.
For $p\in P$ let $Q\restriction p=\{q\in Q:q\leq p\}$ and $Q\uparrow
p=\{q\in Q:q\geq p\}$.
Now $Q\leq_\sigma P$ if and only if for all $p\in P$,
$Q\restriction p$ has a countable cofinal subset and
$Q\uparrow p$ has a countable coinitial subset.
\end{defn}

It is clear that for Boolean algebras $A$ and $B$, if $A$ is a subalgebra
of $B$, then $A\leq_\sigma B$ holds in the Boolean algebraic sense
if and only if  $A\leq_\sigma B$ holds for the partial orders.
Now we can extend the notions $\IDP$, $\SEP$, and $\WFN$ to partial
orders.

\begin{defn} For a partial order $P$, $\SEP(P)$ holds if and only if for all
sufficiently large regular $\chi$ there are cofinally
many $M\in\Mgood{\chi}$ with $P\cap M\leq_\sigma P$.

$\IDP(P)$ holds if and only if for all sufficiently large regular $\chi$ and for all
$M\in\Mgood{\chi}$ with $P\in M$, $P\cap M\leq_\sigma P$.

$\WFN(P)$ holds if and only if for all sufficiently large regular $\chi$ and all
$M\preccurlyeq\mathcal H_\chi$ with $\aleph_1\subseteq M$ and $P\in M$,
$P\cap M\leq_\sigma P$.
\end{defn}

It is easy to check that the extended versions of $\SEP$, $\IDP$, and $\WFN$
agree
with the old ones on Boolean algebras. For $\SEP$ this relies on
Corollary \ref{sepvsidp}, for $\WFN$ on Theorem
\ref{charwfn}.

Let us look at the partial orders
$([\kappa]^{\leq\aleph_0},\subseteq)$.
Assuming the consistency of some large cardinal, it was shown in
\cite{fuso} that $\WFN([\aleph_\omega]^{\leq\aleph_0})$
does not follow from GCH.
However, we have

\begin{thm}
\label{f-s}
{\rm CH} implies
$\IDP([\kappa]^{\leq\aleph_0})$ for every cardinal $\kappa$.
\end{thm}
\begin{proof} Assume CH and let $\chi$ be
sufficiently large and
regular.  
Let $M\in\Mgood{\chi}$
be such that $[\kappa]^{\leq\aleph_0}\in M$. Then we have 
$\kappa\in M$.
By Lemma \ref{CH-M},
$M\cap[\kappa]^{\leq\aleph_0}=[\kappa\cap M]^{\leq\aleph_0}$.

Thus it remains to show that $[\kappa\cap
M]^{\leq\aleph_0}\leq_\sigma[\kappa]^{\leq\aleph_0}$.  Suppose
$x\in[\kappa]^{\leq\aleph_0}$.  If $x\not\subseteq\kappa\cap M$, then
$[\kappa\cap M]^{\leq\aleph_0}\uparrow x=\emptyset$.  Thus,
$[\kappa\cap M]^{\leq\aleph_0}\uparrow x$ is either empty or has a
minimal element, namely $x$.

$[\kappa\cap M]^{\leq\aleph_0}\restriction x$ always has a maximal element,
namely $x\cap M$. This finishes the proof of the theorem.
\end{proof}

This theorem can be regarded as a parallel of Theorem \ref{sepcohen}.
It follows that, assuming the consistency of some large cardinal, it is
consistent that there is a partial order that has the $\IDP$, but not
the $\WFN$.

\section{Complete Boolean algebras satisfying $\SEP$}

Just as the $\WFN$, $\SEP$ is hereditary with respect to
order retracts.  A partial order $P$ is an order retract of a partial
order $Q$ if there are order preserving maps
$e:P\to Q$ and $f:Q\to P$ such that $f\circ e=\id_P$.
If $P$ and $Q$ are Boolean algebras and $e$ and $f$ are
(Boolean) homomorphisms,
then we call $P$ a retract of $Q$.

\begin{lemma}\label{retracts} Let $P$ and $Q$ be partial orders such
that $P$ is an order retract of $Q$.  Then $\SEP(Q)$ implies $\SEP(P)$.
\end{lemma}

\begin{proof} Let $\chi$ be large enough and $M\preccurlyeq\mathcal
H_\chi$.
Suppose $M\cap Q\leq_\sigma Q$ and $P,Q\in M$.
We show $P\cap M\leq_\sigma P$.

Since $M$ knows that $P$ is a retract of $Q$, $M$ contains order preserving maps $e:P\to Q$ and
$f:Q\to P$ such that $f\circ e=\id_P$.
Let $p\in P$.  Since $Q\cap M\leq_\sigma Q$, there is a countable set
$C\subseteq Q\cap M$ such that $C$ is cofinal in $Q\cap M\restriction
e(p)$.

{\em Claim.} $f[C]$ is cofinal in $P\cap M\restriction p$.

Let $q\in P\cap M\restriction p$.  Since $e(q)\in Q\cap M\restriction
e(p)$, there is $c\in C$ such that $e(q)\leq c\leq e(p)$.
Now $q\leq f(c)\leq p$, which proves the claim.

By the same argument, $P\cap M\uparrow p$ has a countable coinitial
subset.  This implies $\SEP(P)$.
\end{proof}

If $P$ is a complete lattice and $P$ embeds into $Q$ via $e$, then there
is a map $f:Q\to P$ with $f\circ e=\id_P$, namely the
map defined by $f(q)=\sup\{p\in P:e(p)\leq q\}$ for all $q\in Q$.
Thus, a complete lattice which embeds into a partial order $Q$
with $\SEP(Q)$ also has the property SEP.

Note that if $A$ is a complete Boolean algebra and $A$ embeds into $B$ in
the Boolean algebraic sense, then $A$ is a retract of $B$ by Sikorski's
extension criterion.
If $A$ is an infinite complete Boolean algebra, then $A$ has a maximal
antichain of size $\aleph_0$ and $\mathcal P(\omega)$ embeds into $A$,
again by Sikorski's extension criterion.  Since $\mathcal P(\omega)$ is
complete, it is a retract of $A$.

It follows that $\SEP(\mathcal P(\omega))$ holds if there is any infinite
complete Boolean algebra $A$ with $\SEP(A)$.  Note that for this it is not
necessary to use Sikorski's criterion;  our statement about complete
lattices suffices.

While $\SEP(\mathcal P(\omega))$ is consistent with (but not a theorem of)
ZFC,
$\SEP(\mathcal P(\omega_1))$ fails.
We prove this in a series of lemmas.
Together with
Lemma \ref{retracts} this will imply that all complete Boolean algebras
$A$ with $\SEP(A)$ satisfy the c.c.c.
It should be pointed out that the
proof of $\neg\WFN(\mathcal P(\omega_1))$ given in \cite{fukoshe}
also works for
$\SEP$.  However, we believe that our argument is simpler.

In the following $\chi$ always denotes a sufficiently large regular cardinal.

\begin{lemma}\label{chain} {\rm a)} $\neg\SEP(\omega_2+1)$

{\rm b)} If $\SEP(A)$ holds, then $A$ does not have a chain of order
type $\omega_2$.
\end{lemma}

\begin{proof}
 For a) note that for each $M\in\Mgood{\chi}$, $\omega_2\cap
M$ is an ordinal of cofinality $\aleph_1$.  It follows that for each
$\alpha\in\omega_2\setminus M$, $(\omega_2+1)\cap
M\restriction\alpha$ has uncountable
cofinality.  In particular, $(\omega_2+1)\cap M\not\leq_\sigma\omega_2+1$.

For b) let $A$ be a Boolean algebra such that $\omega_2$ embeds into
$A$.  Clearly, $\omega_2+1$ also embeds into $A$.  But $\omega_2+1$ is a
complete lattice.  Thus $\omega_2+1$ is an order retract of $A$.  Now
$\neg\SEP(A)$ follows from Lemma \ref{retracts} together with part a).
\end{proof}

\begin{lemma}\label{pomega1cnt} $\neg\SEP(\mathcal
P(\omega_1)/\cnt{\omega_1})$
\end{lemma}

\begin{proof}  By Lemma \ref{chain}, it suffices to show that $\mathcal
P(\omega_1)/\cnt{\omega_1}$ has a chain of ordertype $\omega_2$.

For $f,g\in\omega_1^{\omega_1}$ let $f<^* g$ if and only if
$\{\alpha<\omega_1:f(\alpha)\geq g(\alpha)\}$ is countable.
If $(f_\gamma)_{\gamma<\omega_2}$ is a $<^*$-increasing sequence in
$\omega_1^{\omega_1}$, then
$(\{(\alpha,\beta)\in\omega_1\times\omega_1:\beta\leq
f_\gamma(\alpha)\})_{\gamma<\omega_2}$ gives rise to  a strictly
increasing sequence in $\mathcal
P(\omega_1\times\omega_1)/\cnt{\omega_1\times\omega_1}$ of
order type $\omega_2$.
It is thus sufficient to construct a $<^*$-increasing sequence of order
type $\omega_2$ in
$\omega_1^{\omega_1}$.
But this is easy using the natural diagonalization argument to get an
$<^*$-upper
bound for any set $F\subseteq\omega_1^{\omega_1}$ of size $\leq\aleph_1$.
\end{proof}

\begin{lemma}\label{pomega1} $\neg\SEP(\mathcal P(\omega_1))$
\end{lemma}

\begin{proof} Let $M\in\Mgood{\chi}$.
Suppose $\mathcal P(\omega_1)\cap
M\leq_\sigma\mathcal P(\omega_1)$.  We show that
$$(\mathcal
P(\omega_1)/\cnt{\omega_1})\cap M\leq_\sigma\mathcal
P(\omega_1)/\cnt{\omega_1}.$$
This suffices for the lemma since together with Lemma \ref{pomega1cnt}
it implies that there are not
cofinally many $N\in\Mgood{\chi}$ with $\mathcal
P(\omega_1)\cap
N\leq_\sigma\mathcal P(\omega_1).$
Note that
$(\mathcal P(\omega_1)\cap M)/(\cnt{\omega_1}\cap M)$ is essentially
the same
as $(\mathcal P(\omega_1)/\cnt{\omega_1})\cap M$ since $M$ includes a
cofinal
subset of $\cnt{\omega_1}$, namely the countable ordinals.

Let $a\in\mathcal P(\omega_1)$.  Let $C$ be a cofinal subset of $\mathcal
P(\omega_1)\cap M\restriction a$.  Let $D$ be the set of all classes
modulo $\cnt{\omega_1}$ of elements of $C$ and let $\overline a$ be the
class
of $a$ modulo $\cnt{\omega_1}$.

{\em Claim.}  $D$ is cofinal in $\Big((\mathcal
P(\omega_1)/\cnt{\omega_1})\cap
M\Big)\restriction\overline a$.

Let $b\in\mathcal P(\omega_1)\cap M$ be such that $b\setminus a$ is
countable.  Let $\alpha<\omega_1$ be such that $b\setminus\alpha\subseteq
a$.  Since $\alpha\in M$, $b\setminus\alpha\in M$.  Therefore, there is
$c\in C$ such that $b\setminus\alpha\subseteq c\subseteq a$.
The claim clearly follows from this.  This finishes the proof of the
lemma.
\end{proof}

Let $A$ be a complete Boolean algebra not satisfying the c.c.c.  Since $A$ is
complete, $A$ has an maximal antichain of size $\aleph_1$. This
antichain gives rise to an embedding of the
algebra of finite-cofinite subsets of $\omega_1$ into $A$.  Since $A$ is
complete, this embedding extends to all of $\mathcal P(\omega_1)$ by
Sikorski's
extension criterion.
Since $\mathcal P(\omega_1)$ is complete, it follows that $\mathcal
P(\omega_1)$ is a retract of $A$.  Using Lemma \ref{retracts} this gives

\begin{corollary}\label{ccc} Let $A$ be a complete Boolean algebra.  If
$\SEP(A)$ holds, then $A$ satisfies the c.c.c.
\end{corollary}

\section{An example in ZFC}
\label{example-in-zfc}

In this section we show in ZFC that there is a Boolean algebra $A$ which
satisfies $\SEP(A)$ but not $\IDP(A)$.

\begin{defn}\label{Saleph1} Let
$E^{\aleph_2}_{\aleph_1}=\{\alpha<\omega_2:\cf(\alpha)=\aleph_1\}$.
For all $\alpha\in E^{\aleph_2}_{\aleph_1}$ fix an increasing sequence
$(\delta^\alpha_\beta)_{\beta<\omega_1}$ which is cofinal in $\alpha$ and
consists of successor ordinals.
For $S\subseteq E^{\aleph_2}_{\aleph_1}$ let $A^S$ be the Boolean
algebra defined as follows:
Let $A_0^S=\{0,1\}$.  Suppose $\alpha<\omega_2$ is a limit ordinal and
$A^S_\beta$ has already been defined for all $\beta<\alpha$.  Let
$A^S_\alpha=\bigcup_{\beta<\alpha} A^S_\beta$.  Now suppose $A^S_\alpha$
has been defined and
$\alpha\in\omega_2\setminus S$.  Let
$A^S_{\alpha+1}=A^S_\alpha(x_\alpha)$ where $x_\alpha$ is independent
over $A^S_\alpha$.
Suppose that $A^S_\alpha$ has already been defined and $\alpha\in S$.  Let
$A^S_{\alpha+1}=A^S_\alpha(x_\alpha)$ where
$x_\alpha\not\in A^S_\alpha$,
$A^S_\alpha\restriction
x_\alpha$ is generated by $\{x_{\delta_\beta^\alpha}:\beta<\omega_1\}$ and
$A^S_\alpha\restriction-x_\alpha$ is $\{0\}$.  Finally let
$A^S=\bigcup_{\alpha<\omega_2}A^S_\alpha$.
\end{defn}

Fix a sufficiently large regular cardinal $\chi$.
Recall that for every $M\in\Mgood{\chi}$,
$\omega_1\subseteq M$ and thus, $M\cap\omega_2\in\omega_2$.

\begin{lemma}\label{club} Let $X\subseteq\mathcal H_\chi$ be of size
$\aleph_1$.
Then $C_X=\{M\cap\omega_2:M\in\Mgood{\chi}\wedge X\subseteq M\}$ includes
an  $\aleph_1$-club of $\omega_2$, that is, an unbounded set which is
closed  under limits of subsets of cofinality $\aleph_1$.
\end{lemma}

\begin{proof}
By recursion, we define an increasing  sequence
$(M_\alpha)_{\alpha<\omega_2}$ in
$\Mgood{\chi}$ such that
\begin{itemize}\item[(i)] $X\subseteq M_0$,
\item[(ii)] for all $\alpha<\omega_2$, $\alpha\in M_\alpha$, and
\item[(iii)] if $\alpha<\omega_2$ is a limit ordinal of
cofinality $\aleph_1$, then $M_\alpha=\bigcup_{\beta<\alpha}M_\beta$.
\end{itemize}
This construction can be carried out since $\Mgood{\chi}$ is cofinal in
$[\mathcal H_\chi]^{\aleph_1}$ and closed under unions of chains of
length $\omega_1$.
Let $C=\{\omega_2\cap M_\alpha:\alpha<\omega_2\}$.  Then $C$ is unbounded
in $\omega_2$ by (ii) and closed under limits of subsets of cofinality
$\aleph_1$ by (iii).  By (i), $C\subseteq C_X$.
\end{proof}
From Lemma \ref{club} we get

\begin{lemma}\label{stationary} Let $S$  be
a stationary subset of
$\omega_2$ such that $S\subseteq E^{\aleph_2}_{\aleph_1}$.
Then for cofinally many $M\in\Mgood{\chi}$ we have $\omega_2\cap
M\in S$.
\end{lemma}

\begin{proof}
Let $X\in[\mathcal H_\chi]^{\aleph_1}$.  By Lemma \ref{club},
the set $C_X$ includes an $\aleph_1$-club $C$ of $\omega_2$.  Let 
$\overline C$ be the closure of $C$, i.e., $C$ together with all limit 
points of $C$.   Then $\overline C$ is club in $\omega_2$ and 
$C=\overline C\cap E_{\aleph_1}^{\aleph_2}$.  Since $S$ is stationary 
and a subset of $E_{\aleph_1}^{\aleph_2}$, $\overline C\cap S=C\cap S$ 
is non-empty and thus, there is $M\in\mathcal M_\chi$ such that 
$X\subseteq M$ and $M\cap\omega_2\in S$. 
\end{proof}

Using Lemma \ref{stationary}, we can show that for a suitably chosen set $S$,
the Boolean
algebra $A^S$ constructed above satisfies SEP but not IDP.

\begin{thm}\label{firstexample}  There is a Boolean algebra $A$
with $\SEP(A)$ but not $\IDP(A)$.
\end{thm}

\begin{proof}
Fix two disjoint stationary subsets $S_0$ and $S_1$ of $\omega_2$ with
$S_0\cup S_1=E^{\aleph_2}_{\aleph_1}$.  Let $A=A^{S_1}$.

{\em Claim 1:}  $\SEP(A)$

Let $M\in\Mgood{\chi}$ be such that $\alpha=M\cap\omega_2\in S_0$ and
$A\in M$.
Then $A\cap M=A^{S_1}_\alpha$.  Since $\alpha\not\in S_1$ and by the
construction of $A^{S_1}$, $A^{S_1}_\alpha\leq_\sigma A$.
This proves the claim since there are cofinally many $M\in\Mgood{\chi}$
with $M\cap\omega_2\in S_0$ and $A\in M$ by Lemma \ref{stationary}.

{\em Claim 2:}  $\neg\IDP(A)$

By Lemma \ref{stationary}, there is
$M\in\Mgood{\chi}$ such that $\alpha=M\cap\omega_2\in S_1$ and $A\in M$.
As above, $A\cap M=A^{S_1}_\alpha$.  By the construction of $A^{S_1}$,
$A^{S_1}_\alpha\not\leq_\sigma A$.
In other words, $M$ witnesses the failure of
$\IDP(A)$.
\end{proof}

\section{$\SEP(\mathcal P(\omega))$ does not imply $\IDP(\mathcal
P(\omega))$}\label{section-seven}

In this section we use the idea of the proof of Theorem
\ref{firstexample}
to construct a model of set theory where $\SEP(\mathcal P(\omega))$ holds
while $\IDP(\mathcal P(\omega))$ fails.

\begin{thm} It is consistent that
$\SEP(\mathcal P(\omega))$ holds but $\IDP(\mathcal P(\omega))$ fails.
\end{thm}

\begin{proof}
Suppose the ground model $V$ satisfies CH and
let $S_0$ and $S_1$ be as in the proof of Theorem \ref{firstexample}.
For each $\alpha\in E^{\aleph_2}_{\aleph_1}$ let
$(\delta_\beta^\alpha)_{\beta<\omega_1}$ be as in Definition
\ref{Saleph1}.

Our strategy is to perform a finite support iteration of c.c.c.~forcings
over $V$ of length $\omega_2$ where we add only Cohen reals most of the
time.  However, at
stage $\alpha\in S_1$ we add a new subset $x_\alpha$ of
$\omega$ such that the Cohen reals added at the stages
$\delta_\beta^\alpha$, $\beta<\omega_1$, are almost contained in
$x_\alpha$.
Note that we consider the Cohen reals to be subsets of $\omega$.
This construction should be viewed as the forcing version of the
construction in the proof of Theorem \ref{firstexample}.

We now define the iteration
$(P_\alpha,Q_\alpha)_{\alpha<\omega_2}$.
The underlying sets of the $Q_\alpha$'s will be absolute, but not the
orders.  
Thus we will not define each $Q_\alpha$ as a $P_\alpha$-name
but as the underlying set of $Q_\alpha$ in $V$, also named
$Q_\alpha$, together with a $P_\alpha$-name 
$\dot{\leq}_\alpha$ for the order on
$Q_\alpha$.

For all $\alpha\not\in S_1$ let $Q_\alpha$ be Cohen forcing,
i.e., $\Fn(\omega,2)$.  Let $\dot{\leq}_\alpha$ be the canonical name
for the usual order on $\Fn(\omega,2)$, i.e., reverse inclusion.
For $\alpha\in S_1$ and $\beta<\omega_1$
let $\dot{x}_{\delta_\beta^\alpha}$ be a $P_{\delta^\alpha_\beta+1}$-name
for the Cohen real added by $Q_{\delta^\alpha_\beta}$.
Set $$Q_\alpha=\{(f,F):f\in{}^{<\omega}2,
F\in[\{\delta_\beta^\alpha:\beta\in\omega_1\}]^{<\aleph_0}\}$$ and let
$\dot{\leq}_\alpha$ be a name
for a relation $\leq$ on $Q_\alpha$ such that
$(f,F)\leq(f^\prime,F^\prime)$ if and only if $f^\prime\subseteq f$,
$F^\prime\subseteq F$, and for all
$\delta\in F^\prime$, if $n\in\dot{x}_\delta\cap\dom
(f\setminus f^\prime) $, then $f(n)=1$.

As usual, for each
$\alpha<\omega_2$ let $P_\alpha$ be
the finite support iteration of the $Q_\beta$, $\beta<\alpha$, where each
$Q_\alpha=(Q_\alpha,\dot{\leq}_\alpha)$ is considered as a $P_\alpha$-name
for the appropriate
partial order.  Let $P$ be the direct limit of the $P_\alpha$,
$\alpha<\omega_2$. For convenience, by the absoluteness of the elements of
the $Q_\alpha$'s, we may assume that the elements of each $P_\alpha$ and
of $P$ are elements of
$\prod_{\beta<\omega_2}Q_\beta$ with finite support.  For each condition
$p\in P$ let $\supt(p)$ be its support.

Note that $P$ is c.c.c.~since the $Q_\alpha$'s are $\sigma$-centered.
Let $G$ be $P$-generic over $V$.
For $\alpha<\omega_2$ let $G_\alpha=G\cap P_\alpha$.

The easier part of the proof of the theorem is to show that
in $V[G]$, $\IDP(\mathcal P(\omega))$ fails.
To see this, we need

\begin{lemma}\label{notsigma} For $\alpha\in S_1$, $V[G]\models\mathcal
P(\omega)\cap V[G_\alpha]\not\leq_\sigma\mathcal P(\omega)$.
\end{lemma}

\begin{proof} We argue in $V[G]$.  Let $\alpha\in S_1$.  Note that
$\mathcal P(\omega)\cap V[G_\alpha]\leq_\sigma\mathcal
P(\omega)$ if and only if $(\mathcal P(\omega)\cap
V[G_\alpha])/fin\leq_\sigma\mathcal P(\omega)/fin$.

For $y\in\mathcal P(\omega)$ let $\bar{y}$ be the
equivalence class of $y$ modulo $fin$.
Let $x$ be the subset of $\omega$ generically added by
$Q_\alpha$.

{\em Claim.}
$(\mathcal P(\omega)\cap V[G_\alpha])/fin\restriction\bar x$ is generated
by the classes modulo $fin$ of the Cohen reals
$\{x_{\delta_\beta^\alpha}:\beta<\omega_1\}$ added by the
$Q_{\delta_\beta^\alpha}$'s.  In particular, $(\mathcal
P(\omega)\cap
V[G_\alpha])/fin\not\leq_\sigma\mathcal P(\omega)/fin$.

It follows from the construction of $Q_\alpha$ that for all
$\beta<\omega_1$, $\bar x_{\delta_\beta^\alpha}\leq\bar{x}$.
Let $a\in\mathcal P(\omega)\cap V[G_\alpha]$.
Suppose that $a$ is not almost included in the union of a finite subset
of $\{x_{\delta^\beta_\alpha}:\beta<\omega_1\}$.  Then for every
$n\in\omega$ the set of
conditions in $Q_\alpha$ which force that there is $m\geq n$ such that
$m\in a$ but $m\not\in x$ is easily seen to be dense in $Q_\alpha$.
It follows that $a$ is not almost included in $x$.  This shows the
claim and finishes the proof of Lemma \ref{notsigma}.
\end{proof}

By Lemma \ref{stationary}, in $V$ there is $M\in\Mgood{\chi}$ containing
$P$ such that $\alpha=M\cap\omega_2\in S_1$.
Now $M[G]\cap\omega_2=\alpha$ since $P$ is c.c.c.  If $\dot{x}\in M$ is a
$P$-name for a subset of $\omega$, then $M$ also contains a nice $P$-name
$\dot{y}$ for the same subset of $\omega$.
By c.c.c., $\dot{y}$ only uses
countably many conditions from $P$.  
These conditions are already
contained in $P_\beta$ for some $\beta<\alpha$.  
Since $V$ satisfies CH
and since $P_\beta\in M$ (and thus $P_\beta\subseteq M$), 
$\dot{y}\in M$ by Lemma \ref{CH-M}.  
It follows
that $\mathcal P(\omega)\cap M[G]=\mathcal P(\omega)\cap V[G_\alpha]$.
Therefore, $M[G]$ shows that $\mathcal P(\omega)$ does not
satisfy $\IDP$ in $V[G]$.

To see that $\SEP(\mathcal P(\omega))$ holds in $V[G]$ we need

\begin{lemma}\label{sigma} If $\alpha<\omega_2$ and $\alpha\not\in S_1$,
then $V[G]\models\mathcal P(\omega)\cap V[G_\alpha]\leq_\sigma\mathcal
P(\omega)$.
\end{lemma}

\begin{proof}
By Lemma \ref{cohenreal}, it is
sufficient to show that for $\alpha\in\omega_2\setminus S_1$ every real
in $V[G]\setminus
V[G_\alpha]$ is contained in a Cohen extension of $V[G_\alpha]$.

Let $\alpha\in\omega_2\setminus S_1$.
Let $x\in\mathcal P(\omega)$, but $x\not\in V[G_\alpha]$.
Then there is an $P$-name $\dot{x}\in V$ for $x$.
By c.c.c., we may assume that
$\dot{x}$ uses only countably many conditions from $P$.
Our plan is to find $P_{\dot{x}}\subseteq P$ such that
\begin{enumerate}
\item $P_\alpha\subseteq P_{\dot{x}}$ and $P_\alpha$ is completely
embedded in $P_{\dot{x}}$,
\item $\dot{x}$ is an $P_{\dot{x}}$-name,
\item $P_{\dot{x}}$ is completely embedded in $P$, and
\item the quotient $P_{\dot{x}}:G_\alpha$ is equivalent to
$\Fn(\omega,2)$.
\end{enumerate}

This suffices for the lemma.  For suppose $P_{\dot{x}}$ is as above.
It is not hard to see that (3) implies that $P_{\dot{x}}:G_\alpha$ is
completely embedded in $P:G_\alpha$.  Note that by (1), it is
reasonable to consider $P_{\dot{x}}:G_\alpha$.
By (2), $\dot{x}$ can be regarded as a $P_{\dot{x}}:G_\alpha$-name.
Thus by (4), $x$ is contained in a Cohen extension of $V[G_\alpha]$
and the lemma follows.

It remains to construct $P_{\dot{x}}$ and to show that it has the
required properties.
For a condition $p\in P$ let
$$\overline\supt(p)=\supt(p)\cup\bigcup
\{F:\exists\beta\in\supt(p)\cap S_1\exists
f(p(\beta)=(f,F))\}.$$
Let $$X=\alpha\cup\bigcup\{\supt(p):
\dot{x}\mbox{ uses the condition }p\}$$
and $P_{\dot{x}}=\{p\in P:\overline\supt(p)\subseteq X\}$.

{\em Claim.}  $P_{\dot{x}}$ has the properties (1)--(4).

(1) $P_\alpha\subseteq P_{\dot{x}}$ follows from the definitions.
$P_{\dot{x}}$ can be viewed as finite support iteration
which has $P_\alpha$ as an initial segment.  Thus
$P_\alpha$ is completely embedded in $P_{\dot{x}}$.

(2) It follows from the definitions that $\dot{x}$ is an
$P_{\dot{x}}$-name.

(3) We have to show the following:

\begin{itemize}
\item[(i)] $\forall p,q\in P_{\dot{x}}
(p\perp_{P_{\dot{x}}}q\Rightarrow p\perp_P q)$,
\item[(ii)] $\forall p\in P\exists q\in P_{\dot{x}}\forall r\in
P_{\dot{x}}(r\leq q\Rightarrow r\not\perp_P p)$.
\end{itemize}

For (i) observe that for all $p,q\in P$ with
$p\not\perp_P q$
there is $r\in P$ such that $r\leq p,q$ and
$\overline\supt(r)\subseteq\overline\supt(p)\cup\overline\supt(q)$.
Therefore, if $p,q\in P_{\dot{x}}$ are compatible in $P$, then they
are in $P_{\dot{x}}$.

(ii) is what really requires work.  Let $p\in P$.
Let $q^\prime\in P_{\dot{x}}$
be the condition with support $\supt(p)\cap X$
such that for
all $\beta\in(\supt(p)\cap X)\setminus S_1$,
$q^\prime(\beta)=p(\beta)$ and for all $\beta\in\supt(p)\cap X\cap S_1$,
$q^\prime(\beta)=(f,F\cap X)$ where $f$ and $F$ are such that
$p(\beta)=(f,F)$. $q^\prime$ does not yet work for $q$ in (ii).  We have
to extend it a little.

Let $q$ be the condition with the same support
as $q^\prime$ such that $q(\beta)=q^\prime(\beta)$ for all
$\beta\in\supt(q^\prime)\setminus S_1$.  Now fix
$\beta\in\supt(q^\prime)\cap S_1$.
Let  $f$ and $F$ be such that
$q^\prime(\beta)=(f,F)$
Let $m\in\omega$ be such that for
all $\gamma\in\supt(p)\setminus S_1$, $\dom(p(\gamma))\subseteq m$.
Let $g\in 2^{<\omega}$ be such
that $m\subseteq\dom(g)$, $f\subseteq g$, and
$g(n)=1$ for all $n\in\dom(g)\setminus\dom(f)$.
Now set $q(\beta)=(g,F)$.

{\em Subclaim.} $q$ works for (ii).

Let $r\in P_{\dot{x}}$ be such that $r\leq q$.
We have to construct a common extension $s\in P$ of $p$ and $r$.
As above, we build an approximation $s^\prime$ of $s$ first.
For $\beta\in S_1$ with
$p(\beta)=(f,F)$
and $r(\beta)=(f^\prime,F^\prime)$ let $s^\prime(\beta)=(f\cup
f^\prime,F\cup F^\prime)$.  Note that $f\cup f^\prime$ is a
function since by the definition of $q$ and by $r\leq q$,
we even have $f\subseteq f^\prime$ whenever $\beta\in\supt(r)$.
Note that this definition makes sense if
$\beta\not\in\supt(r)\cap\supt(p)$ since the largest
element of $Q_\beta$ is simply $(\emptyset,\emptyset)$
(for $\beta\in S_1$).

For $\beta\in\omega_2\setminus S_1$ let
$s^\prime(\beta)=p(\beta)\cup r(\beta)$.  Again,
$p(\beta)\cup r(\beta)$ is a function since for $\beta\in\supt(r)$,
$p(\beta)\subseteq r(\beta)$ by $r\leq q$ and the definition of $q$.
It is easy to see that $s^\prime$ extends $r$.

It may happen that $s^\prime\not\leq p$.  However, we can extend
$s^\prime$ to a condition $s\leq p$ by adding some Cohen conditions
(deciding more of the Cohen reals involved).
Let $s(\beta)=s^\prime(\beta)$ for all $\beta\in S_1$.
For $\beta\in S_1$ we have to
make sure
that $s\restriction\beta$ forces $s(\beta)$ to be below
$p(\beta)$.

For all $\beta\in S_1$ let $f_\beta$, $F_\beta$,
$f^\prime_\beta$, and $F^\prime_\beta$ be such that
$s(\beta)=(f_\beta,F_\beta)$ and
$p(\beta)=(f_\beta^\prime,F_\beta^\prime)$.
Then for all $\beta\in S_1$ we have
$f^\prime_\beta\subseteq f_\beta$
and $F_\beta^\prime\subseteq F_\beta$.  For all
$\delta\in F^\prime_\beta$ we want to have
$\delta\in\supt(s)$
and $s(\delta)\forces\forall n\in
\dom(f_\beta)\setminus\dom(f^\prime_\beta)
(f_\beta(n)\geq\dot{x}_\delta(n))$.  This can be accomplished.
Just let $z$
be a sufficiently long finite sequence of zeros and put
$s(\delta)=s^\prime(\delta)\frow z$ for every
$\delta\in\bigcup\{F_\beta:\beta\in\supt(s^\prime)\cap S_1\}$.
Note that there are only finitely
many $\delta$'s to be considered.  For every $\beta\in\omega_2$
for which $s(\beta)$ has not yet been defined let
$s(\beta)=s^\prime(\beta)$.

It is straight forward to check that $s$ is a common extension
of $r$ and $p$.  This completes the proof of the subclaim and thus
shows that $P_{\dot{x}}$ is completely embedded in $P$.

(4) Note that any two elements of $P_{\dot{x}}$ that agree on
$[\alpha,\omega_2)$ are equivalent in $P_{\dot{x}}:G_\alpha$,
i.e., they will be identified in the completion of
$P_{\dot{x}}:G_\alpha$.  But since
$\{\delta^\beta_\gamma:\gamma\in\omega_1\}\cap\alpha$ is countable
for all $\beta\in[\alpha,\omega_2)\cap S_1$ and $X\setminus\alpha$
is countable, there are only countably many possibilities
for $p\restriction[\alpha,\omega_2)$ for $p\in P_{\dot{x}}$.
Therefore, the completion of $P_{\dot{x}}:G_\alpha$
has a countable dense subset.  Since below each element of
$P_{\dot{x}}:G_\alpha$ there are two incompatible elements (in
$P_{\dot{x}}:G_\alpha$),
$P_{\dot{x}}:G_\alpha$ is equivalent to $\Fn(\omega,2)$.
This finishes the proof of the claim and of the lemma.
\end{proof}

Since there are cofinally many $M\in(\Mgood{\chi})^V$ with
$M\cap\omega_2\in S_0$ by Lemma \ref{stationary}, the set
$\{M[G]:M\in\Mgood{\chi}\wedge M\cap\omega_2\in S_0\}$ is cofinal in
$(\Mgood{\chi})^{V[G]}$.  As above, for all $M\in(\Mgood{\chi})^V$,
$\alpha=M[G]\cap\omega_2=M\cap\omega_2$ and $\mathcal P(\omega)\cap
M[G]=\mathcal P(\omega)\cap V[G_\alpha]$.
Now it follows from Lemma \ref{sigma} that
$\SEP(\mathcal P(\omega))$ holds in $V[G]$.  This finishes the proof of
the theorem.
\end{proof}

\section{Variants of SEP}
\label{section-eight}

It is tempting to define a new class of partial orders by replacing 
``cofinally many $M\in{\mathcal H}_\chi$'' in the definition of SEP by 
``stationarily many $M\in{\mathcal H}_\chi$''. However, the class of 
partial orders with this modified notion of SEP coincides with 
the class of partial orders with the original SEP.  
Also, one arrives at the same notion if  ``there are cofinally $M\in{\mathcal H}_\chi$'' is weakened to ``there is $M\in{\mathcal H}_\chi$''.

For a partial order $P$ and a regular cardinal $\chi$ such that 
$P\in\calH_\chi$ let
\[ \calM(P,\chi)=\setof{M\in\calM_\chi}{P\in M\wedge P\cap M\leq_\sigma P}.
\]
\begin{thm}
\label{8.1} Let $P$ be any partial order.  Then for $\kappa=\max(\card{\trcl(P)},\aleph_1)$ the following are equivalent:
\begin{enumerate}\item $\SEP(P)$
\item There is a regular cardinal $\chi>\kappa$ such that $\calM(P,\chi)$ is stationary in $[\mathcal H_{\chi}]^{\aleph_1}$.
\item For every regular cardinal $\chi>\kappa$, $\calM(P,\chi)$ is stationary in $[\mathcal H_{\chi}]^{\aleph_1}$.
\item There is a regular cardinal $\chi>2^\kappa$ such that $\calM(P,\chi)$ is non-empty.
\item For every regular cardinal $\chi>2^\kappa$, $\calM(P,\chi)$ is non-empty.
\end{enumerate}
\end{thm}
Clearly, (1) follows from (3) and implies (4).  The remaining part of
Theorem \ref{8.1} is a special case of the following lemma, which 
does not have to do anything with $\sigma$-embeddings.
For a set $A$, a family $\mathcal F\subseteq\mathcal P(A)$, and a 
regular cardinal $\chi$ with $A,\,\mathcal F\in\mathcal H_\chi$ let  
$$\calM(A,\mathcal F,\chi)=\{M\in\calM_\chi:A,\mathcal F\in M\wedge A\cap M\in\mathcal F\}.$$
For a partial order $P$, $\calM(P,\chi)$ is simply
$\calM(P,\{Q\subseteq P:Q\leq_\sigma P\},\chi)$.  
\begin{lemma}\label{generallemma}
Let $A$ be a set and $\mathcal F\subseteq\mathcal P(A)$.  Then for $\kappa=\max(\card{\trcl(A)},\aleph_1)$ the following are equivalent:
\begin{enumerate}
\item There is a regular cardinal $\chi\geq\kappa^+$ 
such that $\calM(A,\mathcal F,\chi)$ is stationary in $[\mathcal H_{\chi}]^{\aleph_1}$.
\item For every regular cardinal $\chi\geq\kappa^+$, $\calM(A,\mathcal F,\chi)$ is stationary in $[\mathcal H_{\chi}]^{\aleph_1}$.
\item There is a regular cardinal $\chi>2^\kappa$ such that $\calM(A,\mathcal F,\chi)$ is non-empty.
\item For every regular cardinal $\chi>2^\kappa$, $\calM(A,\mathcal F,\chi)$ is non-empty.
\end{enumerate}
\end{lemma}
The reason for considering $\kappa^+$ and $2^\kappa$ in the formulation of this
lemma is that $\kappa^+$ is the least cardinal $\chi>\aleph_1$ with $A\in\mathcal H_\chi$ 
and the size of $\mathcal H_{\kappa^+}$ is $2^\kappa$.
The proof of Lemma \ref{generallemma} uses two arguments: one for stepping up in cardinality and one for stepping down.   We start with decreasing cardinals.
Fix $A$, $\mathcal F$, and $\kappa$  as in Lemma \ref{generallemma}.
\begin{lemma}
\label{8.2}
Let $\chi,\mu>\kappa$ be regular cardinals with 
$2^{<\chi}<\mu$. If $\calM(A,\mathcal F,\mu)$ is non-empty, then $\calM(A,\mathcal F,\chi)$ is stationary in $[\mathcal H_\chi]^{\aleph_1}$.
\end{lemma}
\begin{proof}
Suppose that $\calM(A,\mathcal F,\chi)$ is not stationary in 
$[\mathcal H_\chi]^{\aleph_1}$.  
We may assume that $\chi>\kappa$ is minimal with this property.  Let $M\in\calM(A,\mathcal F,\mu)$.   

Since $\card{\mathcal H_\chi}^{\aleph_1}=2^{<\chi}<\mu$, we have  
$[\mathcal H_\chi]^{\aleph_1},\calM(A,\mathcal F,\chi)\in\mathcal H_\mu$ and $\chi$ is definable in $\mathcal H_\mu$ with the parameters $A$ and $\mathcal F$.
Therefore, $\chi\in M$ and $M$ knows that 
$\calM(A,\mathcal F,\chi)$ is not stationary in $[\mathcal H_\chi]^{\aleph_1}$.

It follows that $M$ contains a club $\mathcal C$ of $[{\mathcal H}_\chi]^{\aleph_1}$
which is disjoint from $\calM(A,\mathcal F,\chi)$.   
By elementarity, $M\cap\mathcal H_\chi\subseteq\bigcup(\mathcal C\cap M)$.
Since $\aleph_1\subseteq M$, $\bigcup(\mathcal C\cap M)\subseteq 
M\cap\mathcal H_\chi$.   Since $M\cap[M]^{\aleph_0}$ is cofinal in $[M]^{\aleph_0}$, $\mathcal C\cap M$ is countably directed.  It follows that $\bigcup(M\cap\mathcal C)$ is the union of an increasing chain of length $\omega_1$ of elements of $\mathcal C$.  Therefore, $M\cap\mathcal H_\chi=\bigcup(\mathcal C\cap M)\in\mathcal C$.  It is easily checked that 
$M\cap\mathcal H_\chi\in\calM_\chi$.  

Since $A\subseteq\mathcal H_\chi$, $A\cap M\cap\mathcal H_\chi=A\cap M\in\mathcal F$.  Thus, $M\cap\mathcal H_\chi\in\mathcal C\cap\calM(A,\mathcal F,\chi)$, contradicting the choice of $\mathcal C$.
\end{proof}

\begin{lemma}
\label{8.3}
For all regular cardinals $\chi, \mu>\kappa$ with $\chi<\mu$, if 
$\calM(A,\mathcal F,\chi)$ is stationary in 
$[{\mathcal H}_\chi]^{\aleph_1}$, then $\calM(A,\mathcal F,\mu)$ is non-empty.
\end{lemma}

\begin{proof} 
We use a refined Skolem hull operator to find $\tilde M\in\calM_\mu$ 
with $A,\mathcal F\in \tilde M$ and $\tilde M\cap A\in\mathcal F$.
Fix a well-ordering  $\sqsubset$ on $\mathcal H_\mu$.
For $\alpha<\omega_1$ let $\sk_\alpha$ denote the Skolem hull operator on 
$\mathcal H_\mu$ with respect to the built-in Skolem functions of the structure
$(\mathcal H_\mu,\in,\sqsubset,A,\mathcal F,\sk_\beta)_{\beta<\alpha}$ where 
$A$ and $\mathcal F$ are considered as constants.
For $X\subseteq\mathcal H_\chi$ let 
$$\sk^*(X)=\bigcup\{\sk_\alpha(Y):
Y\in X\wedge\card{Y}\leq\aleph_0\wedge\alpha<\omega_1\}.$$

{\em Claim.} Let  $X\in [\mathcal H_\mu]^{\aleph_1}$ be such that 
$X\cap [X]^{\leq\aleph_0}$ is cofinal in $[X]^{\leq\aleph_0}$.  Then
\begin{itemize}
\item[(i)] $X\subseteq\sk^*(X)$ and $A,\mathcal F\in\sk^*(X)$,
\item[(ii)] $\card{\sk^*(X)}=\aleph_1$,
\item[(iii)] $[\sk^*(X)]^{\leq\aleph_0}\cap\sk^*(X)$ is cofinal in $[\sk^*(X)]^{\leq\aleph_0}$, and
\item[(iv)] $\sk^*(X)\prec\mathcal H_\mu$.
\end{itemize}

Moreover,
\begin{itemize}
\item[(v)] for all $X,Y\subseteq\mathcal H_\mu$ with $X\subseteq Y$, $\sk^*(X)\subseteq\sk^*(Y)$ and
\item[(vi)] 
if $(X_\alpha)_{\alpha<\delta}$ is an increasing sequence of subsets of $\mathcal H_\mu$ and $X=\bigcup_{\alpha<\delta}X_\alpha$, then 
$\sk^*(X)=\bigcup_{\alpha<\delta}\sk^*(X_\alpha).$
\end{itemize}

For (i) let $x\in X$.  By our assumptions on $X$, there is $Y\in X\cap[X]^{\leq\aleph_0}$ such that  $\{x\}\subseteq Y$.  
Now $x\in Y\subseteq\sk_0(Y)$ and thus $x\in\sk^*(X)$.  This shows $X\subseteq\sk^*(X)$.  $A,\mathcal F\in\sk^*(X)$ since $A,\mathcal F\in\sk_0(Y)$ for every countable element $Y$ of $X$.

Statement (ii) follows from the fact that for every countable set $Y\subseteq\mathcal H_\mu$ and every $\alpha<\omega_1$,
$\sk_\alpha(Y)$ is again countable.  

For (iii) let $Y$ be a countable subset 
of $\sk^*(X)$.  For every $n\in\omega$ fix $\alpha_n<\omega_1$ and a countable set $Y_n\in X$ such that $Y\subseteq\bigcup_{n\in\omega}\sk_{\alpha_n}(Y_n)$.
By our assumptions on $X$, there is $Z\in X\cap[X]^{\leq\aleph_0}$ such that 
$\{Y_n:n\in\omega\}\subseteq Z$.
Let $\beta=\sup_{n\in\omega}(\alpha_n+1)$.   Now for every $n\in\omega$, 
$\sk_{\alpha_n}(Y_n)\in\sk_\beta(Z)$ by the choice of $\sk_\beta$.
Since $\sk_\beta(Z)$ is an elementary submodel of $\mathcal H_\mu$ and since the $\sk_{\alpha_n}(Y_n)$ are countable, 
for every $n\in\omega$ we also have 
$\sk_{\alpha_n}(Y_n)\subseteq\sk_\beta(Z)$.   It follows that $Y\subseteq\sk_\beta(Z)$.  Clearly, $\sk_\beta(Z)$ is a countable subset of $\sk^*(X)$.  We are done with the proof of (iii) if we can show $\sk_\beta(Z)\in\sk^*(X)$.  But this is easy.
Just let $Z^\prime\in X\cap[X]^{\leq	\aleph_0}$ be such that $Z\in Z^\prime$.  
Now $\sk_\beta(Z)\in\sk_{\beta+1}(Z^\prime)\subseteq\sk^*(X)$. 

For (iv) it suffices to show that for all finite subsets $F$ of $\sk^*(X)$ there is an elementary submodel $M$ of $\mathcal H_\mu$ such that $F\subseteq M\subseteq\sk^*(X)$.  
Let $F$ be a finite subset of $\sk^*(X)$.  As before, there are $\beta<\omega_1$ and  $Z\in X\cap[X]^{\leq\aleph_0}$ such that $F\subseteq\sk_\beta(Z)$.
By the definition of $\sk^*$, $\sk_\beta(Z)\subseteq\sk^*(X)$.
And  $\sk_\beta(Z)$ is an elementary submodel of $\mathcal H_\mu$.
 
Statements (v) and (vi) follow immediately from the definition of $\sk^*$.
This finishes the proof of the claim.  
 
Now consider the set 
$$\mathcal C=\{M\in[\mathcal H_\chi]^{\aleph_1}:\sk^*(M)
\cap\mathcal H_\chi=M\}.$$
From the properties of $\sk^*$ it follows that $\mathcal C$ is club in $[\mathcal H_\chi]^{\aleph_1}$.  Since $\mathcal M(X,\mathcal F,\chi)$ is stationary, there is $M\in\mathcal C\cap\mathcal M(X,\mathcal F,\chi)$.  Let $\tilde{M}=\sk^*(M)$.  
By the properties of $\sk^*$, $\tilde M\in\mathcal M_\chi$ and $A,\mathcal F\in\tilde M$.  Moreover, $\tilde M\cap X=M\cap X\in\mathcal F$.  In other words, $\tilde{M}\in\mathcal M(X,\mathcal F,\mu)$.
\end{proof}

\begin{proof}[Proof of Lemma \ref{generallemma}]
We start from (3).
Suppose there is a regular cardinal $\chi>2^\kappa$ such that $\calM(X,\mathcal F,\chi)$ is non-empty.   Then, by Lemma \ref{8.2}, $\calM(X,\mathcal F,\kappa^+)$ is stationary in $[\mathcal H_{\kappa^+}]^{\aleph_1}$.  This implies (1).

Now suppose that (1) holds.  Then there is a regular cardinal $\chi>\kappa$ such that $\calM(X,\mathcal F,\chi)$ is stationary in $[\mathcal H_\chi]^{\aleph_1}$.  
By Lemma \ref{8.3}, there are arbitrarily large regular cardinals $\mu$ such that $\calM(X,\mathcal F,\mu)$ is non-empty.  By Lemma \ref{8.2}, this implies
that $\calM(X,\mathcal F,\mu)$ is stationary for every regular $\mu>\kappa$, i.e., (2) holds.  The implications (2)$\Rightarrow$(4) and (4)$\Rightarrow$(3) are trivial.
\end{proof}

At the moment, we do not know whether $\mathfrak a=\aleph_1$ follows
from $\SEP(\mathcal P(\omega))$. However, we can show that a variant of
$\SEP(\mathcal P(\omega))$ which is called
$\SEP^{+-}(\mathcal P(\omega))$ here (see below) implies
$\mathfrak a=\aleph_1$.

In the following let $\chi$ always denote a regular cardinal.
\begin{defn}
\[ \begin{array}{r@{}l}
		\calM^\sqsubset_\chi=\setof{M}{%
			& M\prec \calH_\chi, \cardof{M}=\aleph_1,
				\mbox{ and there is a well-ordering $\sqsubset$
					on $M$ of order}\\
			&\mbox{type $\omega_1$ such that for every $a\in M$,}
			\sqsubset\cap (M_{\sqsubset a})^2\in M\quad}
   \end{array}
\]\noindent
where $M_{\sqsubset a}=\setof{x\in M}{x\sqsubset a}$.
\end{defn}
\begin{defn}
Let $P$ be a partial order.

(1)
$\SEP^+(P)$ if
$\setof{M\in\calM^\sqsubset_\chi}{P\cap M\leq_\sigma P}$ is cofinal
in $[\calH_\chi]^{\aleph_1}$ for every sufficiently large $\chi$.

(2)
$\SEP^{+-}(P)$ if
$\setof{M\in\calM^\sqsubset_\chi}{P\cap M\leq_\sigma P}$ is non-empty
for a sufficiently large $\chi$.

\end{defn}

For $A$, $\mathcal F$, $\chi$ as in the definition of 
$\calM(A,\mathcal F,\chi)$, let 
$$\calM^\sqsubset(A,\mathcal F,\chi)
=\{M\in\calM^\sqsubset_\chi:A,\mathcal F\in M\wedge A\cap M\in\mathcal F\}.$$
Then it is easy to see that Lemma \ref{8.2} with 
$\calM^\sqsubset(A,\mathcal F,\chi)$ in place of 
$\calM(A,\mathcal F,\chi)$ also holds. As in the proof of Lemma 
\ref{generallemma}, we obtain the following equivalence:
$$
\begin{array}{l}
\phantom{\Leftrightarrow\ \ }\SEP^+(P)\\[\jot]
\Leftrightarrow\ \ 
	\setof{M\in\calM^\sqsubset_\chi}{P\cap M\leq_\sigma P}
		\mbox{ is stationary for every sufficiently large }\chi\\[\jot]
\Leftrightarrow\ \ 
	\setof{M\in\calM^\sqsubset_\chi}{P\cap M\leq_\sigma P}
		\mbox{ is non-empty for every sufficiently large }\chi
\end{array}
$$

\begin{lemma}
\label{x1}$\calM^\sqsubset_\chi\subseteq \calM_\chi$.
\end{lemma}
\begin{proof}
Suppose that $M\in\calM^\sqsubset_\chi$ and $\sqsubset$ is a
well-ordering of $M$ as in the definition of $\calM^\sqsubset_\chi$.
For $X\in[M]^{\aleph_0}$ let $x\in M$ be such that
$X\subseteq M_{\sqsubset x}$. Then
$\cardof{M_{\sqsubset x}}\leq\aleph_0$ and $M_{\sqsubset x}\in M$.
This shows that $[M]^{\aleph_0}\cap M$ is cofinal
in $[M]^{\aleph_0}$. Hence $M\in\calM_\chi$.
\end{proof}
\begin{lemma}
\label{x2}
Suppose that $M\in\calM^\sqsubset_\chi$. Then $M$ is internally
approachable in the sense of {\rm\cite{FoMaShe}}, i.e.,

\begin{itemize}
\item[{\rm ($\ast$)}]
$M$ is the union of a continuously increasing chain of countable
elementary submodels $\seqof{M_\alpha}{\alpha<\omega_1}$ of $M$ such
that $\seqof{M_\beta}{\beta\leq\alpha}\in M_{\alpha+1}$ for all
$\alpha<\omega_1$.
\end{itemize}
\end{lemma}
\begin{proof}
Suppose that $M\in\calM^\sqsubset_\chi$ and $\sqsubset$ is a
well-ordering of $M$ as in the definition of $\calM^\sqsubset_\chi$.

Let $x_\alpha\in M$, $\alpha<\omega_1$ be defined inductively such that
\medskip

(0) $M_{\sqsubset x_\alpha}$ is an elementary submodel of
$\pairof{M,{\in},{\sqsubset}}$;

(1) If $\alpha$ is a limit ordinal then $x_\alpha$ is the limit of
$x_\beta$, $\beta<\alpha$;

(2) If $\alpha$ is a successor, say $\alpha=\beta+1$, then $x_\alpha$
is minimal with respect to $\sqsubset$ such that
$x_\beta\sqsubset x_\alpha$,
${\sqsubset}\cap(M_{\sqsubset x_\beta})^2\in M_{\sqsubset x_\alpha}$ and
(0).
\medskip

\noindent Note that, in (2), the construction is possible since
${\sqsubset}\cap(M_{\sqsubset x_\beta})^2\in M$ by the definition of
$M\in\calM^{\sqsubset}_\chi$. By (0) and since
$\seqof{x_\beta}{\beta\leq\alpha}$ is definable in $M$ with the
parameter ${\sqsubset}\cap(M_{\sqsubset x_\alpha})^2$ we have
$\seqof{x_\beta}{\beta\leq\alpha}\in M_{\sqsubset x_{\alpha+1}}$. It
follows that
$\seqof{M_{\sqsubset x_\beta}}{
	\beta\leq\alpha}\in M_{\sqsubset x_{\alpha+1}}$.
Also $M=\bigcup_{\alpha<\omega_1} M_{\sqsubset x_\alpha}$ since
$\sqsubset$ has order type $\omega_1$. Thus
$\seqof{M_{\sqsubset x_\alpha}}{\alpha<\omega_1}$ is a sequence as
required in ($\ast$).
\end{proof}
The property ($\ast$) in Lemma \ref{x2} almost characterizes elements
of $\calM^\sqsubset_\chi$:
\begin{lemma}
\label{x3}
Let $<^*$ be a well-ordering of $\calH_\chi$ of order type
$\cardof{\calH_\chi}$.  If $\seqof{M_\alpha}{\alpha<\omega_1}$ is a
continuously increasing sequence of countable elementary submodels of
$\pairof{\calH_\chi,{\in},<^*}$ such that
$\seqof{M_\beta}{\beta\leq\alpha}\in M_{\alpha+1}$ for all
$\alpha<\omega_1$, then $M=\bigcup_{\alpha<\omega_1}M_\alpha$ is an
element of $\calM^\sqsubset_\chi$.
\end{lemma}
\begin{proof}
For each $x\in M$, let
$\alpha_x=\min\setof{\alpha<\omega_1}{x\in M_{\alpha+1}}$.
Let $\sqsubset$ be the linear ordering on $M$ defined by
\[ x\sqsubset y\ \Leftrightarrow\ \alpha_x<\alpha_y\ \lor\
	(\alpha_x=\alpha_y\land x<^* y).
\]\noindent
Clearly $\sqsubset$ is a well-ordering on $M$. $\sqsubset$ has order
type $\omega_1$ since every initial segment of $M$ with respect to
$\sqsubset$ is countable and $M$ itself is uncountable.

Let $x\in M$.
We show that ${\sqsubset}\cap(M_{\sqsubset x})^2\in M$.
Let $\alpha^*=\alpha_x+2$. 
In $M_{\alpha^*}$,
${\sqsubset}\cap(M_{\sqsubset x})^2$ is definable from
${<^*}\cap(M_{\alpha_x+1})^2$ and $\seqof{M_\beta}{\beta\leq\alpha_x+1}$.
By elementarity and since ${<^*}\cap(M_{\alpha_x+1})^2$ and
$\seqof{M_\beta}{\beta\leq\alpha_x+1}$ are elements of $M_{\alpha^*}$, it
follows that
${\sqsubset}\cap(M_{\sqsubset x})^2\in M_{\alpha^*}\subseteq M$.
\end{proof}
By Lemma \ref{x2} and Lemma \ref{x3},   there are club many
$M\prec \calH_\chi$ of size $\aleph_1$ such that
$M\in\calM^\sqsubset_\chi$ if and only if $M$ is internally
approachable --- namely those $M$ with $M\prec(\calH_\chi,{\in},<^*)$ for some fixed $<^*$ as above.
\begin{lemma}
\label{x4}
$\calM^\sqsubset_\chi$ is stationary in $[\calH_\chi]^{\aleph_1}$.
\end{lemma}
\begin{proof}
Suppose that $\calC\subseteq[\calH_\chi]^{\aleph_1}$ is closed
unbounded. We show that $\calM^\sqsubset_\chi\cap\calC\not=\emptyset$.

Let $<^*$ be a well-ordering of $\calH_\chi$ of order type
$\cardof{\calH_\chi}$. Let $\seqof{M_\alpha}{\alpha<\omega_1}$ be a
continuously increasing chain of countable elementary submodels of
$\pairof{\calH_\chi,\calC,{\in},{<^*}}$ such that
$\seqof{M_\beta}{\beta\leq\alpha}\in M_{\alpha+1}$ for all
$\alpha<\omega_1$. Let $M=\bigcup_{\alpha<\omega_1}M_\alpha$. Then
$M\in\calM^\sqsubset_\chi$ by Lemma \ref{x3}. Since
$\omega_1\subseteq M$, we have $N\subseteq M$ for all
$N\in\calC\cap M$. By elementarity $\calC\cap M$ is a directed system
and $M=\bigcup(\calC\cap M)$. Since $\cardof{M}=\aleph_1$, it follows
that $M\in\calC$.
\end{proof}
\begin{lemma}
\label{x5}
For a partial order $P$

{\rm(1)} $\IDP(P)$ implies $\SEP^+(P)$;

{\rm(2)} $\SEP^+(P)$ implies $\SEP^{+-}(P)$;

{\rm(3)} $\SEP^{+-}(P)$ implies $\SEP(P)$;

\[\WFN(P)\ \ \Rightarrow\ \ \IDP(P)\ \ \Rightarrow\ \
		\SEP^+(P)\ \ \Rightarrow\ \ \SEP^{+-}(P)\ \ \Rightarrow\ \ \SEP(P)\ 
\]\noindent

\end{lemma}
\begin{proof}
(1) follows from Lemma \ref{x1}. 
(2) is 
clear from definitions. (3) follows from Theorem \ref{8.1}.
\end{proof}
\begin{thm}\label{x6}
Assume $\SEP^{+-}(\psof{\omega})$. Then $\mathfrak a=\aleph_1$.
\end{thm}
\begin{proof}
Let $\chi$ be sufficiently large and $M^*\in\mathcal M^\sqsubset_\chi$ be
such that $\psof{\omega}\cap M^*\leq_\sigma\psof{\omega}$. Since
$\cardof{M^*}=\aleph_1$, it is enough to show that there is a
MAD-family $\subseteq M^*$.

Let $\sqsubset$ be a well-ordering of $M^*$ as in the definition of
$\calM^\sqsubset_\chi$.
By Lemma \ref{x2}, there is an increasing sequence
$\seqof{M_\alpha}{\alpha<\omega_1}$ of countable elementary submodels
of $M^*$ such that\ $\bigcup_{\alpha<\omega_1}M_\alpha=M$ and
$\seqof{M_\beta}{\beta\leq\alpha}\in M_{\alpha+1}$ for all
$\alpha<\omega_1$.

Let $\seqof{a_\alpha}{\alpha<\omega_1}$ be such that
\begin{itemize}
\item[(1)] $\setof{a_n}{n\in\omega}$ is a partition of $\omega$ with
$\seqof{a_n}{n\in\omega}\in M_0$;
\item[(2)] For $\alpha\geq\omega$,
$a_\alpha\in[\omega]^{\aleph_0}\cap M_{\alpha+1}$ is minimal (with
respect
to $\sqsubset$) with the following property:
\\
($\alpha$) $a_\alpha$ is almost disjoint to each
$a_\beta$, $\beta<\alpha$.
\\
($\beta$)
$\forall x\in[\omega]^{\aleph_0}\cap M_\alpha\ \Big(\,
	\forall u\in[\alpha]^{<\aleph_0}\
		(\ \cardof{x\setminus\bigcup_{\beta\in u}a_\beta}=\aleph_0\,)
		\ \rightarrow\ \cardof{a_\alpha\cap x}=\aleph_0\,\Big)$.

\end{itemize}
By (1) and (2), $\seqof{a_\beta}{\beta<\alpha}$ is definable in
$M_{\alpha+1}$ using the parameters ${\sqsubset}\cap(M_\alpha)^2$,
$\seqof{M_\beta}{\beta\leq\alpha}\in M_{\alpha+1}$.
Hence $\seqof{a_\beta}{\beta<\alpha}\in M_{\alpha+1}$
and so the inductive construction of $a_\alpha$ satisfying ($\alpha$) and
($\beta$) can be accomplished in $M_{\alpha+1}$.

By (1) and (2)($\alpha$),
$\setof{a_\beta}{\beta<\omega_1}$ is pairwise almost disjoint. To show
that it is maximal, suppose that it were not. Then there is some
$b\in[\omega]^{\aleph_0}$ such that\ $b$ is almost disjoint to all of
$a_\alpha$'s. Let
$\setof{b_n}{n\in\omega}\subseteq \psof{\omega}\cap M^*$ be a
countable generator of $(\psof{\omega}\cap M^*)\uparrow b$. Let
$\alpha^*<\omega_1$ be such that\
$\setof{b_n}{n\in\omega}\subseteq M_{\alpha^*}$.
Since $a_{\alpha^*}$ and $b$ are almost disjoint there is
$n^*\in\omega$ such that\ $\cardof{b_{n^*}\cap a_{\alpha^*}}<\aleph_0$. By
(2)($\beta$), there is a $u\in[\alpha]^{<\aleph_0}$ such that\
$\cardof{b_{n^*}\setminus\bigcup_{\beta\in u}a_\beta}<\aleph_0$. As
$b\subseteq b_{n^*}$, it
follows that
$\cardof{b\setminus\bigcup_{\beta\in u}a_\beta}<\aleph_0$. But this is
a contradiction  to the choice of $b$.
\end{proof}

Irreversibility of the implications in Lemma \ref{x5} 
cannot be proved in ZFC:
\begin{lemma}
Assume {\rm CH}. 
Then $\calM_\chi=\calM^\sqsubset_\chi$. 
In particular,  for every partial order $P$, 
$\SEP(P)$ if and only if $\SEP^+(P)$.
\end{lemma}
\begin{proof}
$\calM^\sqsubset_\chi\subseteq\calM_\chi$ by Lemma \ref{x1}.
To show $\calM_\chi\subseteq\calM^\sqsubset_\chi$, suppose
$M\in\calM_\chi$. Then by Lemma \ref{CH-M}
$[M]^{\aleph_0}\subseteq M$. Let $\sqsubset$ be an arbitrary
well-ordering of $M$ of order type $\omega_1$. 
Then, for every $x\in M$,
$M_{\sqsubset x}$ and ${\sqsubset}\cap(M_{\sqsubset x})^2$ are
countable subsets of $M$ and hence, by CH, elements of $M$. 
This shows that
$M\in\calM^\sqsubset_\chi$.
\end{proof}
Even under $\neg$CH, $\SEP$ and $\SEP^+$ can be equivalent for partial 
orders with an ``$\aleph_2$-version'' of $\IDP$. 

Let us say that a partial order $P$ has 
the $\aleph_2$-$\IDP$ if for any sufficiently large $\chi$ 
and $M\prec\calH_\chi$, if $\cardof{M}=\aleph_2$, $P\in M$ and 
$[M]^{\aleph_1}\cap M$ is cofinal in $[M]^{\aleph_1}$ then 
$P\cap M\leq_{\aleph_2}P$ where $P\leq_{\aleph_2}Q$ is defined just as 
in Definition 4.1 with ``countable'' there replaced by ``of 
cardinality $<\aleph_2$''. 

Note that every partial order of cardinality $\leq \aleph_2$ has 
the $\aleph_2$-$\IDP$. 
\begin{thm}
\label{x7}
Assume $\Box_{\omega_1}$. For any partial order $P$ 
with the $\aleph_2$-$\IDP$, 
$\SEP(P)$ if and  
only if\/ $\SEP^+(P)$.
\end{thm}
For the proof of Theorem \ref{x7} we use the fact that the class of
partial orders with SEP is closed under $\leq_\sigma$-suborders.
\begin{lemma}
\label{x8}
For partial orders $P$ and $Q$, if $\SEP(P)$ and $Q\leq_\sigma P$ then
$\SEP(Q)$.
\end{lemma}
\begin{proof}
Fix a sufficiently large regular $\chi$. 
It is enough to show that,
for every $M\prec\calH_\chi$ with  $P,Q\in M$, if
$P\cap M\leq_\sigma P$ then $Q\cap M\leq_\sigma Q$.
To see this, let $x_0\in Q$ and we show that
$Q\cap M\restriction x_0$ has a countable cofinal subset. 
(That
$Q\cap M\uparrow x_0$ has a countable coinitial subset can be proved
similarly.) 
By $P\cap M\leq_\sigma P$ there is a countable
set $X\subseteq (P\cap M)\restriction x_0$ \st\ $X$ is cofinal in
$(P\cap M)\restriction x_0$. By $M\models Q\leq_\sigma P$ and
elementarity, for every $x\in X$ we can find  $X_x\in M$  \st\
$M\models\mbox{``\,}X_x\mbox{ is cofinal in }Q\restriction x\mbox{\,''}$.
Then $X_x\subseteq M$ and $X_x$ is a countable cofinal subset of
$(Q\cap M)\restriction x$. 
Letting $Y=\bigcup_{x\in X}X_x$, we have
$Y\subseteq (Q\cap M)\restriction x_0$ and $Y$ is countable.
$Y$ is a cofinal subset
of $(Q\cap M)\restriction x_0$: 
If
$y\in (Q\cap M)\restriction x_0$, then in particular
$y\in (P\cap M)\restriction x_0$. 
Hence there is  $x\in X$ such that $y\leq x$. 
Since
$M\models y\in Q\restriction x$, there is $x'\in X_x\subseteq Y$ \st\
$y\leq x'$.
\end{proof}
Note that the proof above actually shows that each of the variants of
SEP considered above (and  also IDP and WFN) is closed under 
$\leq_\sigma$-suborders.
\begin{proof}[Proof of Theorem \ref{x7}]
If $\cardof{P}<\aleph_2$, then the assertion of the theorem is trivial.
Hence we may assume $\cardof{P}\geq\aleph_2$. 
If $\SEP^+(P)$ then $\SEP(P)$ by
Lemma \ref{x5} (2). 
So we assume $\SEP(P)$ and prove $\SEP^+(P)$.

Let $\chi$ be sufficiently large and let $X$  be an arbitrary element
of $[\calH_\chi]^{\aleph_1}$. 
We show that there is 
$M\in\calM^\sqsubset_\chi$ such that $X\subseteq M$ and
$P\cap M\leq_\sigma P$.

Fix a well-ordering $<^*$ of $\calH_\chi$ of order type
$\cardof{\calH_\chi}$. 
Let
$\calC=\setof{C_\alpha}{\alpha\in Lim(\omega_2)}$ be
a $\Box_{\omega_1}$-sequence.

Let $\seqof{M_\alpha}{\alpha<\omega_2}$ and
$\seqof{a_{\alpha,\gamma}}{\alpha<\omega_2,\,\gamma<\omega_1}$ be
sequences defined inductively so that they satisfy the following
conditions:
\begin{itemize}%
\item[(0)] $\seqof{M_\alpha}{\alpha<\omega_2}$ is a continuously
increasing sequence of elementary submodels of
$\pairof{\calH_\chi,{\in},<^*}$ of cardinality $\aleph_1$.
\item[(1)] $\omega_1$, $X\subseteq M_0$, $P$, $\calC\in M_0$.
\item[(2)]  For all $\alpha<\omega_2$, 
$\seqof{a_{\alpha,\gamma}}{\gamma<\omega_1}$ is an
enumeration of $M_\alpha$.
\item[(3)]  For all $\beta<\omega_2$ we have
$\seqof{M_\alpha}{\alpha\leq\beta},
{<^*}\cap(\bigcup_{\alpha\leq\beta}M_\alpha)^2,
\seqof{a_{\alpha,\gamma}}{\alpha\leq\beta,\,\gamma<\omega_1}
\in M_{\beta+1}$.
\item[(4)]  For all $\alpha<\omega_2$, $P\cap M_{\alpha+1}\leq_\sigma P$.
\end{itemize}
Let $M=\bigcup_{\alpha<\omega_2}M_\alpha$. 
Then $P\cap M\leq_{\aleph_2} P$ by the $\aleph_2$-$\IDP$\ of $P$.
From (4) it follows that $P\cap M\leq_\sigma P$. 
Hence by Lemma 
\ref{x8}, $\SEP(P\cap M)$. 
It follows 
from Lemma \ref{stationary} that there is
$\alpha^*\in E^{\omega_2}_{\omega_1}$ such that
$P\cap M_{\alpha^*}\leq_\sigma P\cap M\leq_\sigma P$. 
Since $X\subseteq M_{\alpha^*}$ by
(1), the proof is complete if we can show the following:

{\em Claim.} $M_{\alpha^*}\in\calM^\sqsubset_\chi$.

Let $C=C_{\alpha^*}$. $C$ is a cofinal subset of $\alpha^*$ of order
type $\omega_1$. 
Let $(\xi_\alpha)_{\alpha<\omega_1}$ be strictly increasing enumeration
of $C$. 
For each limit ordinal $\alpha<\omega_1$
there is $\beta<\alpha^*$ such that $\xi_\alpha\in M_\beta$. 
Since
$C_{\xi_\alpha}=\setof{\xi_\gamma}{\gamma<\alpha}$ by coherence and
since $C_{\xi_\alpha}\in M_\beta$, we have
$\setof{\xi_\gamma}{\gamma<\alpha}\in M_\beta\subseteq M_{\alpha^*}$.
Hence
\begin{enumerate}
\item[($\ast$)] for all $\alpha<\omega_1$,
$\setof{\xi_\gamma}{\gamma<\alpha}\in M_{\alpha^*}$.
\end{enumerate}
\noindent
Let $\varphi:\omega_1\rightarrow\omega_1\times\omega_1$;
$\alpha\mapsto\pairof{\varphi_0(\alpha),\varphi_1(\alpha)}$ be a
surjection such that $\varphi\in M_0$.

We now define a continuously increasing sequence
$\seqof{N_\alpha}{\alpha<\omega_1}$ of countable elementary submodels
of $M_{\alpha^*}$ such that
\begin{itemize}%
\item[(5)] $a_{\xi_{\varphi_0(\alpha)},\varphi_1(\alpha)}$,
$\seqof{N_\beta}{\beta\leq\alpha}\in N_{\alpha+1}$ for all
$\alpha<\omega_1$;
\item[(6)] $N_{\alpha+1}$ is the countable elementary submodel of
$M_{\alpha^*}$ with $N_{\alpha+1}\in M_{\alpha^*}$ which is minimal
with respect to $<^*$ satisfying (5).
\end{itemize}
That this construction is possible can be seen as follows:
By ($\ast$) and since the predicate
``$N_\alpha\prec M_{\alpha^*}$'' can be replaced by
``$N_\alpha\prec M_\eta$'' for a sufficiently large ordinal $\eta<\alpha^*$, each
initial segment of
$\seqof{N_\alpha}{\alpha<\omega_1}$ is definable in $\calH_\chi$ 
with parameters in $M_{\alpha^*}$ and
hence is an element of $M_{\alpha^*}$.

By (5)
$\bigcup_{\alpha<\omega_1}N_\alpha=M_{\alpha^*}$ and
$\seqof{N_\beta}{\beta\leq\alpha}\in N_{\alpha+1}$ for all
$\alpha<\omega_1$.
From Lemma \ref{x3} it follows
that $M_{\alpha^*}\in\calM^\sqsubset_\chi$. 
This finishes the proof of the claim and hence of the theorem.
\end{proof}
\begin{corollary}
Suppose that $\Box_{\omega_1}$ holds and 
$\psof{\omega}$ has the  
$\aleph_2$-$\IDP$ (in particular this is the case if $2^{\aleph_0}=\aleph_2$).
Then $\SEP(\psof{\omega})$ implies $\mathfrak a=\aleph_1$.
\end{corollary}
\begin{proof}
Under the assumptions, if
$\SEP(\psof{\omega})$ then we have
$\SEP^+(\psof{\omega})$ by Theorem \ref{x7}. 
Hence
$\mathfrak{a}=\aleph_1$ by Theorem \ref{x6}.
\end{proof}

\section{Conclusion}
As we have mentioned in the introduction, 
large cardinals are necessary to construct a Boolean algebra
$A$ with the IDP but without the WFN.    
In this sense, IDP and WFN are pretty much the same and it is not surprising
that all the interesting set theoretic consequences of $\WFN(\psof{\omega})$
that have been discovered so far already follow from $\IDP(\psof{\omega})$.
Looking at the proofs of the known consequences of $\WFN(\psof{\omega})$ or
$\IDP(\psof{\omega})$,
it turns out that most of the time $\SEP(\psof{\omega})$ is enough
to derive these consequences.  
An exception could be the equality $\mathfrak a=\aleph_1$, which is not known to
follow from $\SEP(\psof{\omega})$, but which follows from
$\IDP(\psof{\omega})$.  
The natural open question is whether $\SEP(\psof{\omega})$ $+$ 
$\mathfrak a>\aleph_1$
is consistent.  

One nice feature of $\SEP(\psof{\omega})$ is that it holds in Cohen models.  
This does not have to be true for $\WFN(\psof{\omega})$ (assuming large cardinals).
We do not know about $\IDP(\psof{\omega})$.  
As it turns out, $\SEP(\psof{\omega})$ is relatively robust under slight changes of the
definition.
It does not matter whether we demand the existence of a single elementary 
submodel of $\calH_\chi$ with certain properties, or of stationarily many, or of cofinally many.  
Therefore it is interesting to know that the strongest variant of SEP along these lines, IDP, is strictly stronger than SEP.  
In some sense, we get the best possible result here. 
There is  (in ZFC) a Boolean algebra with SEP but without IDP and   
it is consistent that $\psof{\omega}$ itself is an example.

%
%

\end{document}